\documentclass[12pt, twoside, epsf]{article}

\usepackage{color,graphicx,times}
\usepackage{amsmath, amssymb, graphics}

\newcommand{\mathsym}[1]{{}}
\newcommand{\unicode}[1]{{}}

\makeatletter
\long\def\M#1{\leavevmode\setbox\@tempboxa\hbox{#1}\@tempdima\fboxrule
    \advance\@tempdima \fboxsep \advance\@tempdima \dp\@tempboxa
   \hbox{\lower \@tempdima\hbox
  {\vbox{\hrule \@height \fboxrule
          \hbox{  \hskip\fboxsep
          \vbox{\vskip\fboxsep \box\@tempboxa\vskip\fboxsep}\hskip
                 \fboxsep\vrule \@width \fboxrule}%
                  }}}}
\makeatother

\let \ttorg \tt \def \tt{\ttorg \obeyspaces}

\begin{document}

\newcommand{\Across}{\raisebox{-0.25\height}{\includegraphics[width=0.5cm]{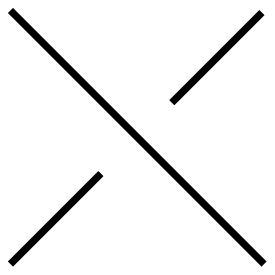}}}
\newcommand{\Asmooth}{\raisebox{-0.25\height}{\includegraphics[width=0.5cm]{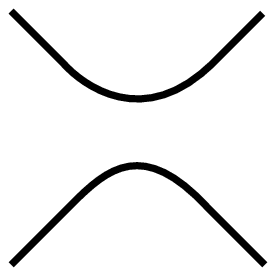}}}
\newcommand{\Bsmooth}{\raisebox{-0.25\height}{\includegraphics[width=0.5cm]{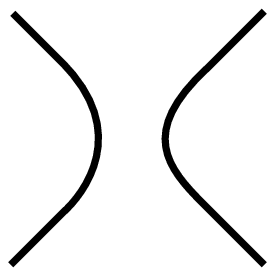}}}
\newcommand{\Rcurl}{\raisebox{-0.25\height}{\includegraphics[width=0.5cm]{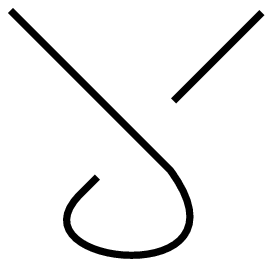}}}
\newcommand{\Lcurl}{\raisebox{-0.25\height}{\includegraphics[width=0.5cm]{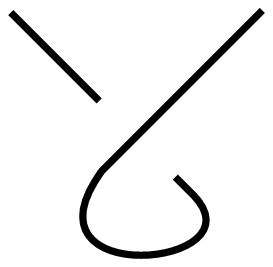}}}
\newcommand{\Arc}{\raisebox{-0.25\height}{\includegraphics[width=0.5cm]{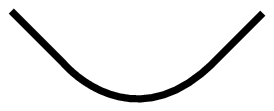}}}

\date{}

\title{\bf Simplicial Homotopy Theory, Link Homology and Khovanov Homology}

\author{Louis H. Kauffman \\
 Department of Mathematics, Statistics and Computer Science \\
 University of Illinois at Chicago \\
 851 South Morgan Street\\
 Chicago, IL, 60607-7045}

\maketitle
  
\thispagestyle{empty}

\section{Introduction}
The purpose of this note is to point out that simplicial methods and  the well-known Dold-Kan construction in 
simplicial homotopy theory \cite{Dold,DoldPuppe,May,Curtis,Lamotke,Weibel} can be fruitfully applied to convert link homology theories into homotopy theories. Dold and Kan \cite{Dold} prove that there is a functor $\Gamma$ from the  category $\cal{C}$ of chain complexes over a commutative ring $R$ with unit to the category $\cal{S}$ of simplicial objects over $R$ such that chain homotopic maps in $\cal{C}$ go to {\it homotopic} maps in $\cal{S}.$ Furthermore, this is an equivalence of categories
$$\Gamma: \cal{C} \longrightarrow \cal{S}.$$
 In this way, given a link homology theory, we construct a mapping  $$T: \cal{L} \longrightarrow \cal{S}$$ taking
link diagrams $\cal{L}$ to a category of simplicial objects $\cal{S}$ such that up to looping or delooping,
link diagrams related by Reidemeister moves will give rise to homotopy equivalent simplicial objects, and the homotopy groups of these objects will be equal to the link homology groups of the original link homology theory. The construction is independent of the particular link homology theory, applying equally well to Khovanov Homology and to Knot Floer Homology and other theories of these types. \bigbreak

A simplifying point in producing a homotopy simplicial object in relation to a chain complex occurs when the chain complex is itself derived (via face maps) from a simplicial object that satisfies the Kan extension condition.
See Section 3 of the present paper. Under these circumstances one can use that simplicial object rather than apply the Dold-Kan functor to the chain complex. We will give examples of this situation in regard to Khovanov homology. The results that we describe here are similar, for Khovanov homology, to the results of \cite{ET,ELST,LS1,LS2}.  
We will investigate the exact relationship with these results in separate papers. And we will investigate detailed working out of this correspondence in separate papers. The purpose of this note is to announce the basic relationships for using simplicial methods in this domain. Thus we do more than just quote the Dold-Kan Theorem.
We give a review of simplicial theory and we point to specific constructions, particularly in relation to Khovanov homology, that can be used to make simplicial homotopy types directly.
\bigbreak

The main ideas in this paper are summarized in the figures. Figure~\ref{prism}  and Figure~\ref{horns} illustrate how degeneracy operators are essential in forming cartesian products and homotopy extension for simplicial objects. Figure~\ref{Figure 2} Illustrates the Khovanov Category for the trefoil knot. The Khovanov Category associates a category to a knot or link diagram. As we shall see, Khovanov homology is essentially the homology 
of a functorial image of this category as a Frobenius module category. Figure~\ref{morphism}  shows how a sequence of $n$ composable morphisms in any category has the structure of an $n$ dimensional simplex. This is
the key notion behind the homology of small categories and the notion of the nerve of a category (the simplicial structure on all composable sequences of morphisms). Figure~\ref{simplexcat} illustrates how the barycentric 
subdivision of the simplex associated with a simplex category corresponds to the nerve of that category. Figure~\ref{Figure 3} illustrates the cube category that stands behind the Khovanov category of Figure~\ref{Figure 2}.
Figure~\ref{cubeplex} and Figure~\ref{cubecat} show how an $n$-cube category fits inside an $n$-simplex category.\\

Section 2 is a review of homology and chain homotopy.
Section 3 is a concise review of simplicial theory.
Section 4 describes the Dold-Kan correspondence.
Section 5 describes the application of the Dold-Kan Theorem to link homology theories.
Section 6 describes the simplicial background to the homology and cohomology of small categories. We show how to construct, for any small category $C$, a simplicial object $NerveC$ defined in terms of composable sequences 
of morphisms in the category. A functor $F$ from $C$ to a module category gives rise to a simplical module $ZF(NerveC).$ This simplicial module satisfies the Kan extension condition and is therefore a homotopy object. The homology of the category $C$ with coefficents in $F$ is equal to the homology of this simplicial object: $H_{*}(C,F) \cong H_{*}(F(NerveC)).$ Furthermore, as observed by Jozef Przytycki and Jin Wang \cite{Jozef1,Jozef2,Wang},
for a simplex category $Cat[n]$, its barycentric subdivision is identical with $NerveCat[n].$ We make this correspondence precise in the Subdivision Theorem in Section 6. The upshot of this result is that one can identify 
the homology of a simplex category with coefficients in $F$ with the homology of a homotopy simplicial module object: $H_{*}(Cat[n], F) \cong H_{*}(F(NerveCat[n]))$. In Section 7 we show how Khovanov homology for a knot or link $K$, $Kho_{*}(K)$, can be identified in the form   $H_{*}(Cat[n], F)$ and hence one can associate a homotopy object $ZF(NerveCat[n])$ for the Khovanov homology of the link $K.$ Khovanov homology is the homology of a small category (with appropriate coefficients) and  thereby the homology of a simplicial object. This method is direct and circumvents 
using the Dold-Kan Theorem to construct the homotopy object. We give this construction in Section 7 at the end of the section, preceded by a discussion of the embedding of the cube categories in simplex categories and a discussion of another direct construction at the level of the simplex category.\\ 

Section 8 (an appendix) of the paper describes the relationship of Khovanov homology to the bracket polynomial and to the states of the bracket. We include this appendix to provide background for the original construction of Khovanov homology. Here is a quick description of how Khovanov homology relates to the present paper. Khovanov \cite{Kho} associates a small category to a knot diagram. The method of this association is illustrated in Figure~\ref{Figure 2} where we show the bracket states for the trefoil knot arranged with arrows indicating a single re-smoothing of a state. Smoothings of type $A$ have an arrow to smoothings of type $B$ providing the directionality for the category. The objects of this category are the bracket states themselves, but the pattern behind them is a cube category, as shown in Figure~\ref{Figure 3}. One wants to measure the Khovanov category in order to obtain topological information about the knot or link that generates it. Khovanov constructs a functor to a module category by associating a module to each bracket state and module maps to the re-smoothing
arrows. The modules are chosen so that a homology calculation related to the category is invariant (with certain grading conventions) under isotopies of the knots and links. We discuss how to replace the category by a 
simplicial space whose homotopy type carries this homology. Cubical structures can be embedded in simplicial structures as in Figure~\ref{cubeplex}. Simplicial structures are directly related to composable sequences of 
morphisms in any category. See Figure~\ref{morphism}. In this way simplicial topology is a background structure in the category. Sequences of morphisms in the module category form the generating elements for
a simplicial object in the technical sense (with both face and degeneracy operators -- See Section 3) that can be taken as the appropriate homotopy object behind the Khovanov homology. This is the object $ZF(NerveCat[n])$
referred to above.\\

 In this paper we do not discuss the construction of a spectrum of homotopy types, but this is needed to compare categories for links that are related by Reidemeister moves. We describe in Section 3 the appropriate looping and delooping functors on simplicial objects and this can be used to compare the chain homotopies between Khovanov complexes of links that are related by Reidemeister moves. See \cite{K,K1,BN1,BN2} for a description of these chain homotopies (they involve grading shifts that correspond here to the looping and delooping operations on simplicial objects). It is our intent in a further paper to compare
chain homotopies of the Khovanov complex with homotopies of the corresponding simplicial objects. This paper is intended to be a first paper in a series of papers on this topic. For that reason, it contains what the author regards as the necessary background for building further work in these directions.\\

We end this introduction with a question: Are there smaller and more accessible spatial models than the simplicial models indicated here? In particular, let $K$ be a link diagram. Let $\Lambda(K)$ denote the nerve of the image of the Khovanov category for
$K$ under the Frobenius algebra functor that associates a module to each object of the Khovanov category for $K.$ The homology of $\Lambda(K)$ can be taken as the Khovanov homology of $K.$ The simplicial object
$\Lambda(K)$ is simpler than the free abelian group object $Z\Lambda(K)$ whose homotopy type is a product of Eilenberg-Maclane spaces for the homology of $\Lambda(K).$ In this paper, we have promoted 
$Z\Lambda(K)$ as a possibly useful homotopy type for Khovanov homology of $K.$ But $\Lambda(K)$ already carries this homology. It is a very interesting question to ask about the homotopy type of the realization
$|\Lambda(K)|.$ This is a space that carries the Khovanov homology. What are its homotopy groups? How do they behave under isotopies of the link $K?$ We are indebted to Stephan Klaus for this line of questioning.\\

\noindent {\bf Acknowledgement.} We thank A. K. Bousfield for an early conversation, Chris Gomes for many conversations and for his participation in the Quantum Topology Seminar at UIC, Jonathan Schneider, Jozef Przytycki and John Bryden for many conversations about this project. We thank Stephan Klaus of the Mathematisches Forschungsinstitut Oberwolfach for helpful conversations and we thank the Mathematisches Forschungsinstitut Oberwolfach for their hospitality while much of this work was completed.\\
 
 \section{The Category of Chain Complexes}
In this section we review the category of chain complexes and the concept of chain homotopy.
\bigbreak

Let $C$ denote a chain complex
over a commutative ring $R$ with unit.
Thus there is a module $C_{n}$ over $R$  for each $n = 0,1,2,\cdots $ and maps  of modules
$$\partial_{k}: C_{k}\longrightarrow C_{k-1}$$ for $k >0.$ There is also an augmentation mapping\
$$\epsilon: C_{0} \longrightarrow R.$$ The augmentation is a map of modules where we regard $R$ as a module over itself by left multiplication.  It is assumed that $\partial_{k} \circ \partial_{k+1} = 0$ for 
all $k >0$ and that $\epsilon \circ \partial_{1} = 0.$ The {\it homology groups} of $C$ are the groups defined by the equation $$H_{k}(C) = Kernel(\partial_{k})/Image(\partial_{k+1})$$ for $k>0$.
We define $H_{0}(C) = Kernel(\epsilon)/Image(\partial_{1}).$  An elment $x \in C_{k}$ is said to have 
{\it degree $k$.}
\bigbreak

A map of chain complexes $f:C\longrightarrow C'$ is a collection of module maps $f_{n} : C_{n} \longrightarrow C'_{n}$ such that  $\partial'_{n} \circ f_{n} = f_{n-1} \circ \partial_{n}$ for all $n>1$ and 
 $\epsilon'  \circ f_{0} = id_{R}\circ \epsilon$ where $id_{R}$ is the identity map on the ring $R.$  A map of chain complexes induces a homomorphism of the corresponding homology groups: $f_{*}: H_{n}(C) \longrightarrow H_{n}(C').$
 \bigbreak

Let $I$ denote the (unit interval) chain complex with $I_{0} = \langle e_{0}, e_{1} \rangle $ and 
$I_{1} = \langle e \rangle$ where $\langle \cdots \rangle$ means the module generated by the contents of the brackets. We take $\partial_{1}(e)= e_{0} - e_{1}$ and $\epsilon(e_{0}) = \epsilon(e_{1}) = 1.$
\bigbreak 

Recall that two chain maps $f_{0}, f_{1}: C \longrightarrow C'$ are said to be {\it chain homotopic}  if there is 
a chain map $D: I \otimes C \longrightarrow C'$ such that $D(e_{0} \otimes x) = f_{0}(x)$ and 
$D(e_{1} \otimes x) = f_{1}(x)$ for all $x \in C$ where $x \in C_k$ for some $k.$ The boundary map for the tensor product of chain complexes is given by the formula $\partial(a \otimes b) = \partial(a) \otimes b + (-1)^{|a|} a \otimes \partial(b)$ where $|a|$ denotes the degree of $a.$ It is easy to see that chain homotopic maps induce identical maps on homology. One says that two chain complexes $C$ and $C'$ are  {\it chain homotopy equivalent} if there are chain maps $f: C\longrightarrow C'$ and 
$g:C'\longrightarrow C$ such that $f \circ g$ and $g \circ f$ are each chain homotopic to the identity map of their respective domains. It follows that chain homotopy equivalent complexes have isomorphic homology.
 \bigbreak
 
 Let $\cal{C}$ denote the category of chain complexes over the ring $R$ as we have discussed them above. The objects in this category are the chain complexes themselves. The morphisms in the category are the chain maps We have the notion of chain homotopy of maps in this category and the notion of chain homotopy equivalence of objects in the category.
 \bigbreak

  \section{Recalling Simplicial Theory}
  We begin by recalling the structure of abstract simplices. An abstract (non-degenerate) $n$-simplex is denoted by  $$\langle i_1\, i_2 \, \cdots \, i_n \rangle$$
 where $i_1 < i_2 < \cdots < i_n$ and the $i_k$ are elements of an ordered index set, possibly infinite. For purposes of illustration and without loss of generality, we will use 
 $\langle 0 \,1\, 2 \,\cdots \, n \rangle$ as a representative abstract simplex.\\
 
 Face operators $d_i , i = 0, \cdots n, $ applied to a non-degenerate $n$-simplex, produce $n$ distinct $n-1$ simplices. The face operator $d_i $ removes the $i$-th entry of the given simplex.
 Thus $$d_i \langle 0 \,1\, 2 \,\cdots \, (i-1) \, i \, (i+1) \, \cdots \, n \rangle  = \langle 0 \,1\, 2 \,\cdots \, (i-1) \,\, (i+1) \, \cdots \, n \rangle$$ and we often write 
 $$d_i  \langle 0 \,1\, 2 \,\cdots \, n \rangle = \langle 0 \,1\, 2 \,\cdots \,  \hat{i} \, \cdots \, n \rangle$$ where $\hat{i}$ denotes the elimination of the entry at the $i$-th place.\\
 
 We generalize non-degenerate abstract simplices by relaxing the condition $i_1 < i_2 < \cdots  < i_n$ to $i_1 \le i_2  \le \cdots  \le i_n$ so that the sequence in the simplex is monotone and can have equal 
 adjacent elements. Then one has {\it degeneracy operators} $s_i $ that, when applied to an $n$-simplex (non-degenerate or degenerate) produce an $n+1$-simplex. The operator $s_i $ repeats the 
 $i$-th entry. Thus $$s_i \langle 0 \,1\, 2 \,\cdots \,( i-1) \,\, i \,\, (i+1) \, \cdots \, n \rangle = \langle 0 \,1\, 2 \,\cdots \,( i-1) \,\, i \,\, i \,\, \, (i+1) \, \cdots \, n \rangle.$$
 It is then not hard to see that the following identities hold, where $I$ denotes the identity mapping.\\
 
 \noindent{\bf Simplicial Identities}
 \begin{enumerate}
 \item $d_i d_j = d_{j-1} d_i $ for $i < j.$
 \item $s_i s_j = s_{j+1} s_{i}$ for $i \le j.$
 \item $d_i s_j = s_{j-1} d_i $ for $i < j.$
 \item $d_i s_j = Id$ for $i=j, j+1.$
 \item $d_i s_j = s_j d_{i-1}$ for $i > j+1.$
 \end{enumerate}
 
 \noindent Given a category $C,$ a {\it simplicial object} $X$ over $C$ consists in a series of objects $$X_{0}, X_{1}, \cdots. X_{n}, \cdots  $$ in $C$ and morphisms 
 \begin{enumerate}
 \item $d_{i}: X_{n} \longrightarrow X_{n-1}$ for $0 \le i \le n,$
 \item $s_{i}: X_{n} \longrightarrow X_{n+1}$ for $0 \le i \le n$
 \end{enumerate}
 satisfying the simplicial identities indicated above. (An object is said to be {\it semi-simplicial} if it satisfies the identities involving only face operators and there are no degeneracy operators.)
 If $C$ is a category of sets, then we say that a simplicial object over $C$ is a {\it simplicial set}.\\
 
 If $X$ and $Y$ are simplicial objects over a category $C$ then a {\it map of simplicial objects} $$f: X \longrightarrow Y$$ consists in a collection of maps
 $$f_{n}: X_{n} \longrightarrow Y_{n}$$ for each $n = 0 ,1, 2, \cdots $ such that these maps all commute with the face and degeneracy maps for $X$ and for $Y.$\\

 \noindent {\bf Chain Complex and Homology.} If $X$ is a simplicial set, and $R$ is a commutative ring with unit, then we let $C_{i}(X)$ be the free module over $R$ generated by the elements of $X_{i}.$ We extend the mappings $d_{i}$ and $s_{i}$ linearly
 so that we have maps $d_{i}: C_{i}(X) \longrightarrow C_{i-1}(X)$ and $s_{i}:C_{i}(X) \longrightarrow C_{i+1}(X).$ It is not hard to see that this makes the collection $ C_{*}(X)= \{ C_{i}(X) | i \ge 0\} $ into a new simplicial set
 whose objects are modules over $R.$ Furthermore we can define $$\partial_{n} : C_{n}(X) \longrightarrow C_{n-1}(X)$$ by the formula $$\partial_{n} = \sum_{k=0}^{n} (-1)^{k} d_{k}.$$ It is not hard to see that 
 $$\partial_{n-1} \partial_{n} = 0$$ so that $C_{*}(X)$ is a chain complex. The homology of this complex will be written as $H_{*}(X; R).$\\
 
 \noindent{\bf Singular Complex of a Space.} Let $\Delta_{n}$ denote the standard geometric $n$-simplex. $$\Delta_{n} = \{(t_{0},t_{1},\cdots, t_{n}) | 0 \le t_{i} \le 1, \sum_{i=0}^{n} t_{i} = 1 \}.$$
 Define $\delta_{i}: \Delta_{n-1} \longrightarrow \Delta_{n}$ and $\sigma_{i}: \Delta_{n+1} \longrightarrow \Delta_{n}$ by
 $$\delta_{i}((t_{0},\cdots t_{n-1})) = (t_{0},\cdots, t_{i-1}, 0 , t_{i}, \cdots , t_{n-1}),$$
 $$\sigma_{i}((t_{0},\cdots t_{n+1})) = (t_{0},\cdots, t_{i}+  t_{i+1}, \cdots , t_{n+1}).$$
 These maps are dual to the face and degeneracy maps that we have discussed for simplicial structure. The {\it singular complex $S(X)$} of a space $X$ is the simplicial set defined via
 $S(X)_{n}$ equals the set of continuous maps of the geometric $n$- simplex $\Delta_{n}$ to $X.$ If $f:\Delta_{n} \longrightarrow X,$ then we define $$d_{i}f = f \circ \delta_{i}$$ and 
 $$s_{i}f = f \circ \sigma_{i}.$$ This gives $S(X)$ the structure of a simplicial set. The homology we have defined above for $S(X)$,  $H_{*}(S(X); R)$ is the standard singular homology of the space $X.$\\
 
\noindent{\bf Geometric Realization.} If $K$ is a simplicial set, then there is a {\it geometric realization of $K$} denoted $|K|.$ We recall the construction \cite{May}.
Give $K$ the discrete topology and form the disjoint union $$\bar{K} = \sqcup_{n \ge 0} (K_{n} \times \Delta_{n}).$$ Define an equivalence relation $\sim$ on $\bar{K}$ via 
$$(d_{i}k, t) \sim (k, \delta_{i}(t))$$ and $$(s_{i}k, t) \sim (k, \sigma_{i}(t)).$$ Here $\delta_{i}$ and $\sigma_{i}$ are defined as above in the description of the singular complex.
The geometric realization $|K|$ is the quotient of this disjoint union by the equivalence relation $\sim:$  $$|K| = \bar{K}/\sim.$$ Thus we associate a geometric simplex to each element of the simplicial set, 
and these geometric simplices are glued together in accordance with the face and degeneracy relations in the simplicial set. It is then not hard to see that $H_{*}(|K|,R)$ is isomorphic with 
$H_{*}(K;R),$ as we have defined it above. In fact $|K|$ is a $CW$ complex and behaves well with respect to simplicial homotopy, which we shall define below.\\
\\
   
 \noindent{\bf Simplex Examples.} A good example of a simplicial set (simplicial object over the category of sets) is the {\it $n$-simplex complex} $\Delta[n].$ The abstract simplex $P_{n} = \langle 0 \, 1\, \cdots \, n \rangle$
 is a member of $\Delta[n]_{n}.$ All faces and all degeneracies generated from $P_{n}$ constitute the whole of $\Delta[n].$ Note that this means that all simplices in $\Delta[n]_{k}$ for $k > n$ are degenerate.
 $\Delta[n]_{0}$ consists in $ \{ \langle 0 \rangle, \langle 1 \rangle,
 \cdots , \langle n \rangle \}.$ All other degrees contain many degenerate abstract simplices. \\
 
 For example $\Delta[0]_{0} = \{ \langle 0 \rangle \},$ $\Delta[0]_{1} = \{ \langle 0 \, 0 \rangle \},$
 $\Delta[0]_{2} = \{ \langle 0 \, 0 \, 0 \rangle \}, $ and so on, with  $\Delta[0]_{k}$ consisting (for $k>0$) in a degenerate abstract simplex with $k+1$ zeroes.\\
 
 For a second example, consider the first three degrees of $\Delta[1].$ All simplices in $\Delta[1]$ are degenerate above degree $1.$
 $$\Delta[1]_{0} = \{ \langle 0 \rangle , \langle 1 \rangle \},$$
  $$\Delta[1]_{1} = \{ \langle 0 \, 0 \rangle,  \langle 0 \, 1 \rangle,  \langle 1 \, 1 \rangle \},$$ 
  $$\Delta[1]_{2} = \{ \langle 0 \, 0 \, 0 \rangle, \langle 0 \, 0 \, 1 \rangle, \langle 0 \, 1 \, 1 \rangle, \langle 1 \, 1 \, 1 \rangle   \}.$$ \\
  We will return to the $n$-simplices after describing products.\\

  \noindent{\bf Products of Simplicial Objects}\\
  Given two simplicial objects $X$ and $Y$ over a category $C$ with products,  their {\it cartesian product} is the object $X \times Y$ defined by taking 
 $$(X \times Y)_{n} = X_{n} \times Y_{n},$$
 and 
 $$d_{i}(x,y) = (d_{i}(x), d_{i}(y)),$$
 $$s_{i}(x,y) = (s_{i}(x), s_{i}(y)).$$
 (We give the formulas in set-theoretic form, but they naturally generalize to an arbitrary category with products.)\\

\begin{figure}
     \begin{center}
     \begin{tabular}{c}
     \includegraphics[width=6cm]{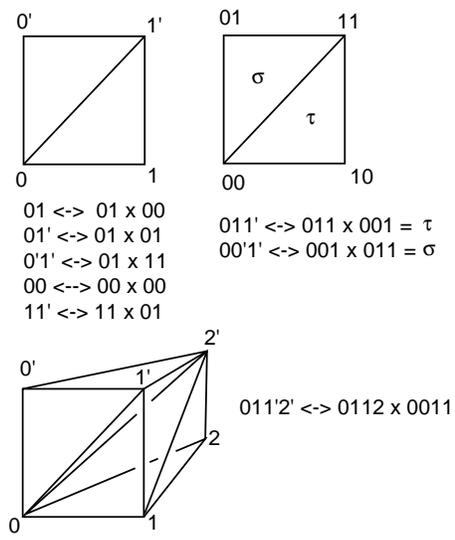}
     \end{tabular}
     \caption{\bf Prism Construction}
     \label{prism}
\end{center}
\end{figure}

The role of the degeneracy operators is crucial for the construction of products.  A key example is $\Delta[n] \times I$ where $I$ denotes $\Delta[1].$
The basic simplices of $\Delta[n] \times I$ consist in the followng with $0 \le a_{0} \le a_{1} \le \cdots \le a_{q} \le n:$
$$( \langle a_{0} \, a_{1} \, \cdots \, a_{q}  \rangle \, , \, \langle 0 \, 0 \, \cdots \, 0 \, 1 \, 1 \, \cdots \, 1 \rangle ).$$
This is abbreviated as 
$$\langle a_{0} \, a_{1} \, \cdots \, a_{i} \, a'_{i+1} \, a'_{i+2} \, \cdots \, a'_{q} \rangle.$$
Here there are $i$ zeroes in the expression $ \langle 0 \, 0 \, \cdots \, 0 \, 1 \, 1 \, \cdots \, 1 \rangle.$
The non-degenerate simplices where  $0 \le a_{0} < a_{1} < \cdots < a_{q} \le n$ correspond to the prismatic decomposition of a standard cartesian product of a geometric $n$-simplex with a unit interval.
This is illustrated in Figure~\ref{prism}. Note how the prismatic decomposition is encoded in the simplicial framework via the degenerate second factors 
 $ \langle 0 \, 0 \, \cdots \, 0 \, 1 \, 1 \, \cdots \, 1 \rangle$ that tell, with $0$ and $1,$ when to take the lower and upper vertices in the prism. Examine Figure~\ref{prism}. \\
 
 \noindent{\bf Homotopy of Simplicial Maps}\\
 If $f:X \longrightarrow Y$ and $g:X \longrightarrow Y$ are two simplicial maps, then we say that $f$ and $g$ are {\it homotopic} if there is a simplicial map $F: X \times I \longrightarrow Y$ such that 
 $F|_{X \times (0)} = f$ and $F|_{X \times (1)} = g.$ In this way the concept of homotopy is carried to simplicial objects. {\it Note that here $(0)$ refers to the simplicial object generated by $\langle 0 \rangle$, and 
 $(1)$ refers to the simplicial object generated by $\langle 1 \rangle.$}\\
 
\begin{figure}
     \begin{center}
     \begin{tabular}{c}
     \includegraphics[width=6cm]{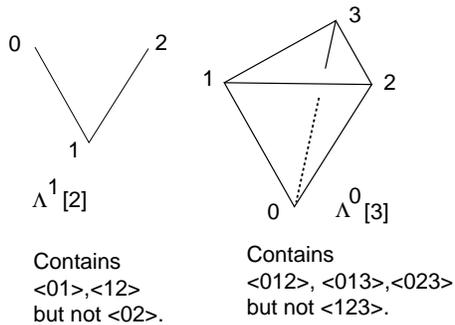}
     \end{tabular}
     \caption{\bf Horns $\Lambda^{k}[n].$}
     \label{horns}
\end{center}
\end{figure}

 One obtains a good homotopy theory when the simplicial objects satisfy the {\it Kan extension condition} which we shall now define.
 Let  the {\it $k$- horn}  $\Lambda^{k}[n]$ denote the subcomplex of $\Delta[n]$ generated by all $d_{i}\langle 0 \, 1 \, \cdots \, n\rangle$ for $i \ne k.$ 
 See Figure~\ref{horns} for an illustration of horns for $n = 2,3.$ A simplicial object is said to {\it satisfy the Kan extension condition} if every map
 $$f:  \Lambda^{k}[n] \longrightarrow K$$ extends to a map $$\hat{f}: \Delta[n] \longrightarrow K.$$ \\
 
 It is easy to verify that the singular complex $S(X)$ of a space $X$ satisfies the Kan extension condition.\\
 
 \noindent{\bf Homotopy Groups.} For a simplicial set $K$ that satisfies the Kan extension condition, one defines the {\it homotopy groups} $\pi_{n}(K,\phi)$  as follows \cite{May}. Let $\phi \in K_{0}$ generate a subcomplex of 
 $K$ having exactly one degenerate simplex in each dimension greater than zero. We let $\phi$ denote either this subcomplex or any of its simplices. 
 Let $\tilde{K_{n}}$ denote the set of  $x \in K_{n}$ such that $d_{i}x = \phi$ for all $ 0 \le i \le n.$ Such an $x$ corresponds to a map $f:\Delta[n] \longrightarrow K$ such that all the faces of $P_{n} = \langle 0 \, 1\, \cdots \, n \rangle$ are taken to $\phi.$ $\pi_{n}(K,\phi)$ consists in the homotopy classes of such maps relative to the subcomplex $\phi,$ and can be seen as $\tilde{K_{n}}/\sim$ where $\sim$ is an appropriate notion of homotopy 
 of simplices in $K$ (See \cite{May}). With this method the definition of the homotopy groups is given in terms of the combinatorial structure of the simplicial set $K.$  In the following we shall often abbreviate 
 $\pi_{n}(K,\phi)$ to $\pi_{n}(K)$ with the subcomplex $\phi$ taken for granted.\\

\noindent{\bf Simplicial Modules and Simplicial Groups.}  A simplicial object $X$ over a {\it module category} (the objects are modules and the morphisms are maps of modules over some given commutative ring) will be called a {\it simplicial module} or 
 a {\it simplicial abelian group} or a {\it simplicial group} (here it is understood that the groups are abelian).\\
 
 \noindent {\bf Theorem \cite{May}.} Simplicial groups satisfy the Kan extension condition.\\
  
  Note that if we have a simplicial set $K,$ then we can form the free abelian simplicial set $ZK$ generated by the simplices of $K.$  The chain complex $C_{*}(K)$ described above where $C_{n}(K)$ is the free module, over the ring of integers $Z,$ on the $n$-simplices of $K$ is obtained from $ZK$ by defining the boundary operator on $ZK$ to be the alternating sum of the face maps from $K$ and extended linearly to $ZK.$
  This makes $ZK = C_{*}(K)$ a simplicial group and hence it satisfies the Kan extension condition, and so can be regarded as an ingredient in a homotopy theory. Furthermore, we have the geometric 
  realization   $|ZK|= |C_{*}(K)|,$ and it follows from the theory of simplicial groups that simplicial homotopies of $C_*(K)$ correspond to homotopies of the realization $|C_{*}(K)|$ as a $CW$ complex.
  It is the case that the homology of the geometric realization $H_{*}(|C_{*}(K)|) = H_{*}(|ZK|) = H_{*}(|K|)$ is isomorphic with the homology of $C_{*}(K)$ as a chain complex defined by its simplicial structure. If it is only homology that we are concerned with then we can take the geometric realization $|K|$ of $K$ itself. It appears to be an open problem to understand the homotopy type of $|K|$, but it is known \cite{May} that 
  $|ZK|$ has the homotopy type of a product of Eilenberg-MacLane spaces of type $K(H_{n}(K),n).$ \\
  
  In the geometric realization of $C_{*}(K)$ as a simplicial object we create a geometric simplex for every element of $C_{*}(K),$ not just for the elements of $K$ but also for all linear combinations of them with 
  coefficients in the ring $R.$ Thus we need to see that given such a linear combination, it is homologous to ``itself" as a single simplicial element. An example may clarify the issue. Suppose that $S$ is a simplicial group
  and that $a$ and $b$ are elements of $S_{1}.$ Let $[a+b]$ denote $(a+b)$ regarded as a single element of $S_{1}$.  Let $$\gamma = s_{0}a + s_{1}b.$$ Then 
  $$d_{0} \gamma = d_{0}s_{0}a + d_{0}s_{1}b = a + s_{0}d_{0}b,$$
  $$d_{1}\gamma = d_{1}s_{0}a + d_{1}s_{1}b = a + b = [a+b],$$
  $$d_{2}\gamma = d_{2}s_{0}a + d_{2}s_{1}b = s_{0}d_{1}a + b.$$
  Thus $$\partial_{2}\gamma = a - [a+b]  + b +(s_{0}d_{0}b + s_{0}d_{1}a),$$ from which it follows that 
  $$\partial_{2}(\gamma  - s_{0}s_{0} d_{0}b - s_{0}s_{0}d_{1}a) = a + b - [a + b].$$
  Thus $[a+b]$ is homologous to $a + b.$ This is actually an instance of the Kan extension condition. 
  We emphasize this point because it may seem unintuitive to make a geometric realization related to a given simplicial set $K$ by realizing the simplicial group $C_{*}(K),$ but the advantage is that 
  $C_{*}(K)$ is in a homotopy category, and the homology of $K$ is preserved even though many new simplices are added in the passage to the simplicial group.\\

 \noindent{\bf Classifying Space.} Let $G$ be a topological group, and define the simplicial set $BG$ with $$BG_{0} = {1}, BG_{1} = G, BG_{2} = G \times G, \cdots, BG_{n} = G^{n}, \cdots .$$ Face and degeneracy operators are defined as follows:
 \begin{enumerate}
 \item $s_{i}(g_{1}, \cdots, g_{n}) = (g_{1}, \cdots, g_{i}, 1, g_{i+1}, \cdots , g_{n}),$
\item  $d_{0}(g_{1}, \cdots, g_{n}) = (g_{2}, \cdots, g_{n})$ 
 \item $d_{i}(g_{1}, \cdots, g_{n}) = (g_{1}, \cdots, g_{i}g_{i+1}, \cdots g_{n}),$ for $0 < i < n,$
 \item $d_{n}(g_{1}, \cdots, g_{n}) = (g_{1}, \cdots , g_{n-1}).$
 \end{enumerate}
 The geometric realization $|BG|$ is called the {\it classifying space} of $G.$ The set $$[X, |BG|]$$ of homotopy classes of maps from a space $X$ to the classifying space $|BG|$ classifies fiber bundles over $X$ with structure group $G.$ One proves that $$\pi_{k}(G) = \pi_{k+1}(BG) =  \pi_{k+1}(|BG|).$$ This classifying space construction can be generalized to the case where $G$ is a simplicial group. See \cite{May}, Chapter 4.\\
 
 \noindent{\bf Path Space and Loop Space.} If $A$ is a simplicial object, the {\it path space} $PA$  \cite{Weibel,May} is the simplicial object defined by $(PA)_{n} = A_{n+1}.$ Here $d_{i}$ on $PA$ is $d_{i+1}$ on $A$ and 
 $s_{i}$ on $PA$ is $s_{i+1}$ on $A.$ The maps $d_{0}: A_{n+1} \longrightarrow A_{n}$ form a simplicial map $PA \longrightarrow A.$ One can show that $PA$ is homotopy equivalent to the constant object
 $A_{0}.$  If $G$ is a group, we write $EG = P(BG)$ where $BG$ is the simplicial classfiying space for $G$ defined above. One shows that $d_{0}:EG \longrightarrow BG$ is a principal fibration (note that $EG$ is 
 contractible) and so the long homotopy sequence of this fibration shows that the homotopy groups of $BG$ are those of $G$ with an index shift. $\pi_{k}(BG) = \pi_{k-1}(G).$\\
 
 Let $A$ be a simplicial object in an abelian category $\cal{A}.$ Let $\Lambda A$ be the simplicial object of $\cal{A}$ that is the kernel of the mapping $d_{0}: PA \longrightarrow A.$ $\Lambda A$ is a loop space for $A.$
 We have $\pi_{k}(\Lambda A) = \pi_{k+1}(A)$ when $A$ is a Kan complex.
 (See \cite{Weibel}). By taking the loop space or the classifying space we can shift the index of the homotopy groups of a given Kan complex accordingly, and this will be of use to us in making constructions later in this paper.
 We shall refer to these (functorial) constructions as {\it looping} and {\it delooping.}
 \\

\begin{figure}
     \begin{center}
     \begin{tabular}{c}
     \includegraphics[width=7cm]{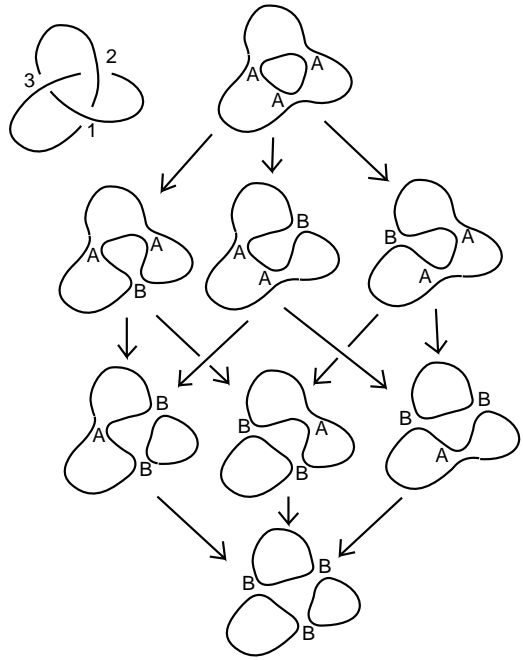}
     \end{tabular}
     \caption{\bf Bracket States and Khovanov Category}
    \label{Figure 2}
\end{center}
\end{figure}

\begin{figure}
     \begin{center}
     \begin{tabular}{c}
     \includegraphics[width=7cm]{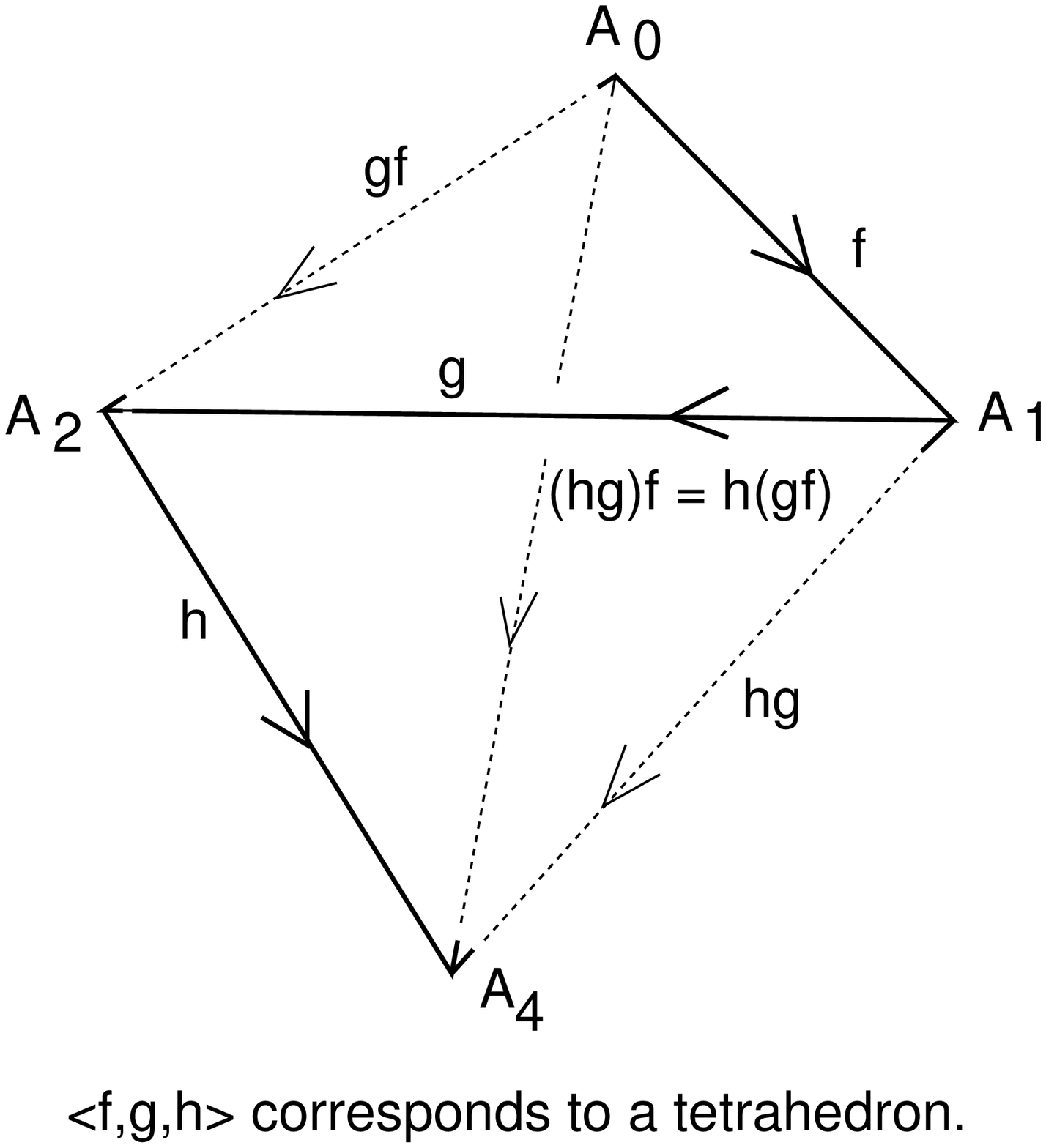}
     \end{tabular}
     \caption{\bf Simplicial Structure of a Morphism Sequence}
    \label{morphism}
\end{center}
\end{figure}

 \begin{figure}
     \begin{center}
     \begin{tabular}{c}
     \includegraphics[width=6cm]{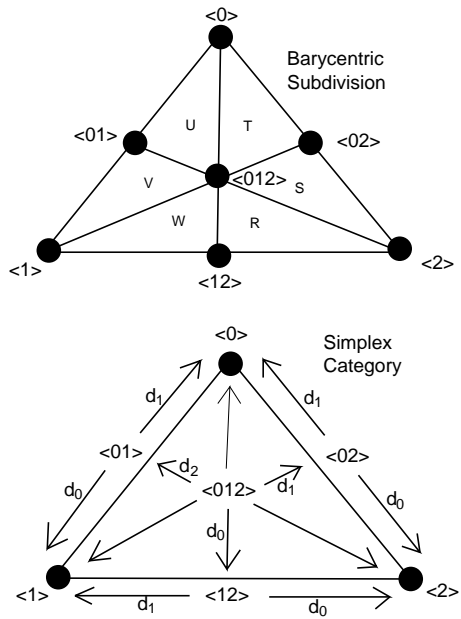}
     \end{tabular}
     \caption{\bf Barycentric Subdivision and the Simplex Category}
     \label{simplexcat}
\end{center}
\end{figure}

\begin{figure}
     \begin{center}
     \begin{tabular}{c}
     \includegraphics[width=7cm]{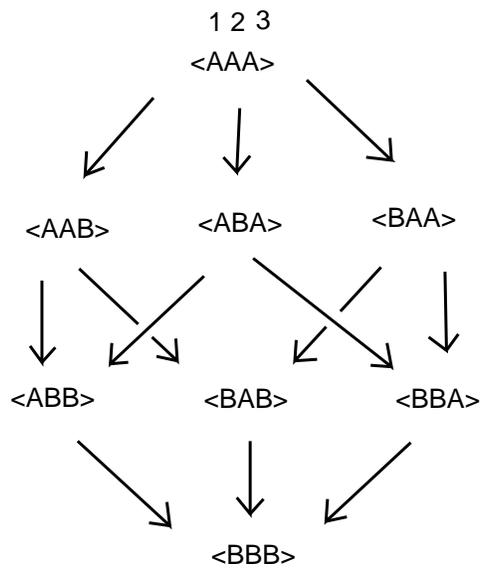}
     \end{tabular}
     \caption{\bf Cube Category}
    \label{Figure 3}
\end{center}
\end{figure}

 \begin{figure}
     \begin{center}
     \begin{tabular}{c}
     \includegraphics[width=6cm]{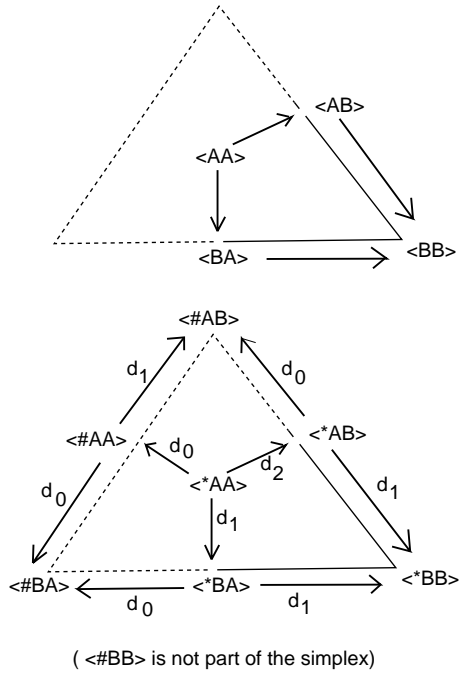}
     \end{tabular}
     \caption{\bf Cube Category Inside Simplex Category}
     \label{cubeplex}
\end{center}
\end{figure}

 \begin{figure}
     \begin{center}
     \begin{tabular}{c}
     \includegraphics[width=6cm]{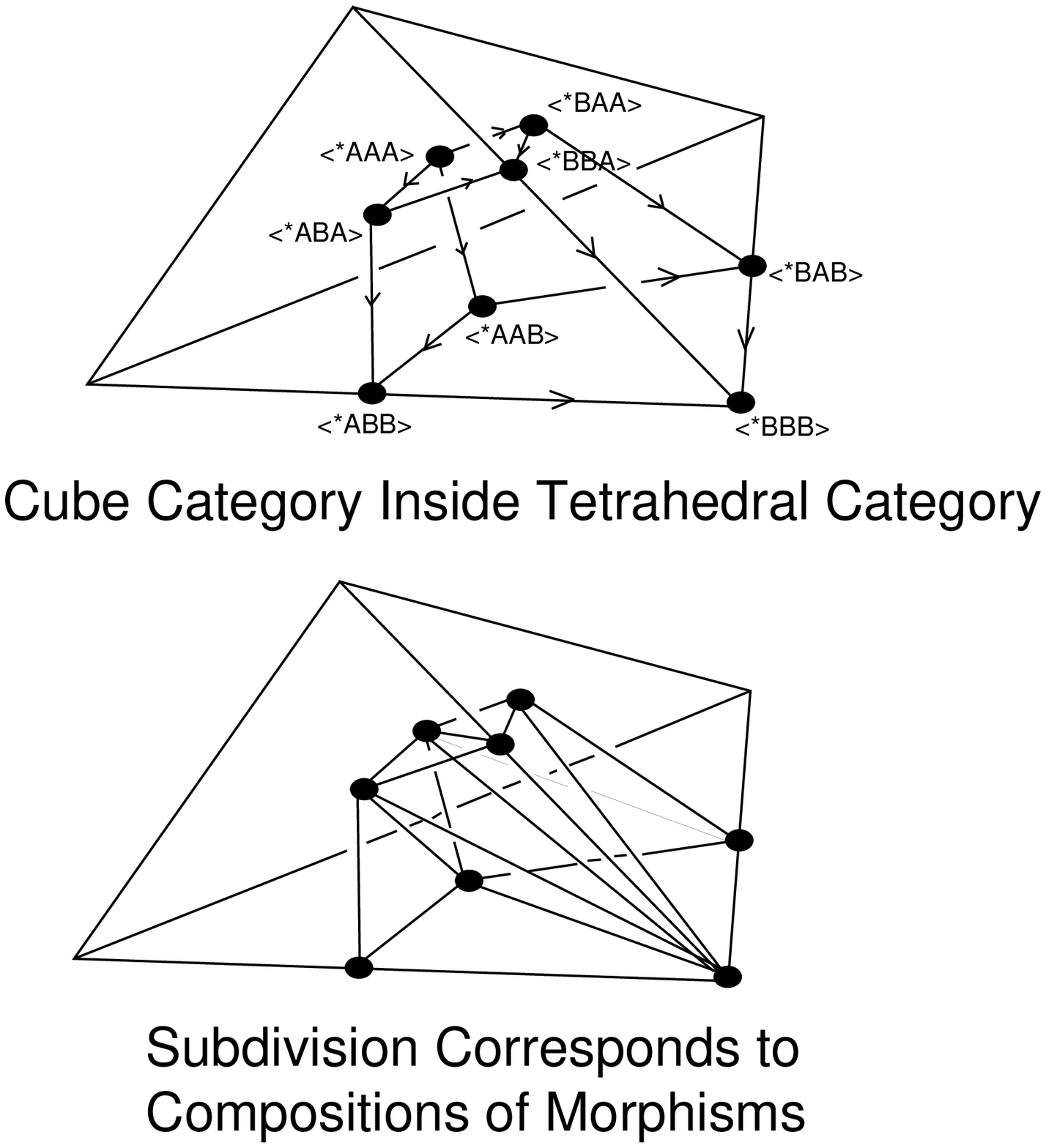}
     \end{tabular}
     \caption{\bf Barycentric Subdivision and the Simplex Category in Dimension Three}
     \label{cubecat}
\end{center}
\end{figure}

 \section{The Dold-Kan Theorem}
 Dold and Kan \cite{Dold} prove that there is a functor $\Gamma$ from the  category $\cal{C}$ of chain complexes over a commutative ring $R$ with unit to the category $\cal{S}$ of simplicial objects over $R$ (the levels are modules over $R$) such that chain homotopic maps in $\cal{C}$ go to {\it homotopic} maps in $\cal{S}.$ Furthermore, this is an equivalence of categories. The inverting functor from $\cal{S}$ to $\cal{C}$ is $N,$ the functor that associates to a simplicial object $S$ its normalized complex (Moore complex)  $N(S)$ where 
 $$N(S)_{n} = \sqcap_{k<n} Kernel(d_{k})$$ where $d_{k}$ is the $k$-th face operator for $S.$
The normalization functor makes any simplicial object over a ring $R$ (all its levels are modules over $R$) into a chain complex since the last face operator at each level can now be seen as a boundary operator. Again, $N$ takes homotopies to chain homotopies. 
\bigbreak

We will describe the functor $\Gamma: \cal{C} \longrightarrow \cal{S}$ below, but first a few more remarks about the functor $N.$ Homotopy groups are defined for simplicial objects over $R$ and it is a fact that 
$$H_{*}(S)  = H_{*}(N(S)) = \pi_{*}(N(S)).$$ Thus the homotopy groups of the Moore complex of a simplicial object over $R$ are the same as the homotopy groups of the orginal simplicial object, and these are the same as its  homology groups.
In the category of simplical objects over $R$, homotopy groups and homology groups coincide.
This means that the Dold-Kan correspondence $\Gamma$ has the property that 
$$\pi_{*}(\Gamma(C)) = H_{*}(N(\Gamma(C)) = H_{*}(C).$$ \\

To avoid excessive notation, we use the symbol $\Delta [n]$ to denote an $n$-simplex object over $R$. This means that its levels are the free module over $R$ generated by the simplices of the original set theoretic object
$\Delta [n]$ that we have previouusly discussed. Then let $$\bar{\Delta}[n] =N(\Delta[n])$$ be the Moore complex for this simplicial object. This is our 
functorial substitute for an $n$-simplex in the category $\cal{S}.$
\bigbreak

Now we can define the functor $\Gamma: \cal{C} \longrightarrow \cal{S}.$ Define
$$\Gamma(C)_{n} = [ \bar{\Delta}[n] \longrightarrow C]$$ where $[ A \longrightarrow B ]$ denotes all chain maps from $A$ to $B$. The face and degeneracy maps for $\Gamma(C)$ are natural and we omit their description here.In this way the Dold-Kan construction associates to a chain complex $C$ a simplicial object $$\Gamma(C)_{*} = [ \bar{\Delta}[*] \longrightarrow C]$$ with the property \cite{Dold} that chain homotopic maps of complexes are taken by $K$ to homotopic maps of the corresponding simplicial objects, and the homology of the chain complexes corresponds functorially to the homotopy groups of the simplicial objects.
\bigbreak

The use of $ \bar{\Delta}[n] $, the Moore complex for the simplicial object $ \Delta [n],$ was the original choice of Dold and Kan \cite{Dold}. It is conceptually simpler but combinatorially more complicated to 
replace   $ \bar{\Delta}[n]$ by $[\Delta[n]] $ where this denotes the chain complex directly associated with $\Delta[n]$ whose boundary operators are the alternating sums of the face operators $d_{k}.$
Thus one can define $\Gamma(C)_{n} = [[\Delta[n]] \longrightarrow C],$ as is done in \cite{Lamotke}. Other choices of definition for the Dold-Kan Functor are available as explained in \cite{DoldPuppe,May,Curtis,Weibel}.
In the course of our further work, we will compare these definitions and how they work in the sense that the functor $K$ induces homotopies of simplicial objects associated with chain homotopies for chain complexes.\\

The basic result of this Dold-Kan correspondence is that there is an equivalence of categories between the category of chain complexes and the category of simplicial abelian groups. At the level of these categories, if $S$ is a simplicial abelian group, then $\Gamma(N(S))$ is a simplicial abelian group isomorphic with $S,$ and if $C$ is a chain complex, then $N(\Gamma(C))$ is a chain complex isomorphic with $C.$
On top of this, the chain homotopy class of a given chain complex $C$ is equivalent to the simplicial homotopy class of its image $\Gamma(C)$ under the Dold-Kan functor. Thus the homotopy type of 
$\Gamma(C)$ classifies the homology type of $C,$ where by {\it homology type} we mean the chain-homotopy type of the complex. \\

As we shall see, in the case of Khovanov homology, there is a natural simplicial abelian group from which the chain complex for Khovanov homology is derived and it is the homotopy type of this simplicial abelian group that becomes a homotopy type for Khovanov homology via the Dold-Kan correspondence.\\

\section{Applications of the Dold-Kan Theorem to Link Homology}
We can be brief. The point of this note is that if there is a homology theory for links with  corresponding  chain complexes $C(L)$ for a link $L$  (for example Khovanov homology), then one can apply the Dold-Kan functor  $\Gamma$ and obtain simplicial objects (or their spatial realizations) such that $\Gamma(C(L))$ has a homotopy theory corresponding to the chain homotopy theory of $C(L).$ In the case of Khovanov homology and other link homology theories, if  one has that links $L$ and $L'$ that are related by Reidemeister moves then the  complexes $C(L)$ and $C(L')$ that are chain homotopy equivalent {\it up to degree shifts  of these complexes}. From this and our discussion of the Dold-Kan Theorem, it follows that $\Gamma(C(L))$ and $\Gamma(C(L'))$ are of the same homotopy type {\it up to the application of looping or delooping}.
Here looping and delooping refer to taking either the loop space construct or the classifying space construct. Taking care of the book-keeping of looping and delooping would associate a spectrum in the sense of homotopy theory  to the link $L.$ The homotopy type of this spectrum would be an invariant of the link. We can say that, up to grading shifts, the homotopy type of $\Gamma(C(L))$ is an invariant of the link. The homotopy groups of $\Gamma(C(L)$ are the same as the Khovanov homology groups (with appropriate attention to grading). It may seem unsatisfactory to resort to the looping and delooping operations, but the possible advantage of the direct use of the Dold-Kan functor is that one can, in principle, investigate the relationship between chain-homotopy on the link homology complex and corresponding homotopies of the simplicial spaces. We shall not attempt to construct spectra in this paper, as there are numerous technical issues in that endeavor (see \cite{Marc}) but here just concentrate, in the next section, on the construction of a simplicial space associated to the Khovanov homology associated with a given link diagram.
\bigbreak

The advantage, from our point of view, in using the Dold-Kan construction  to study homotopy for link homology is that it applies to all such theories  \cite{Kho,BN1,BN2,S,K,K1}. In the case of existing theories there are many specifics to explore and in the case of yet to be formulated theories, we can say that there is the beginning of a homotopy theory of this sort waiting for them. \\

\section{Simplicial Objects, Homotopy and Homology of Categories}
Let $C$ be a small category (i.e. the objects and morphisms form sets).
We now define a simplicial object $NerveC$ associated with the category $C$ as follows.
An element of the $n$-th grading, $NerveC_{n},$ is a sequence $$ \gamma = \langle f_0,f_1,\cdots, f_{n-1} \rangle$$ of morphisms in  $C$ such that 
there are objects $A_0 , A_1 , \cdots , A_n $ so that  
$$f_{k}:A_{k} \longrightarrow A_{k+1}$$ for $k = 0,\cdots, n-1.$\\

\noindent We define face and degeneracy operators by the following formulas.
 \begin{enumerate}
\item $d_{k}\gamma =  \langle f_{0}, \cdots, f_{k-2}, f_{k}\circ f_{k-1},f_{k+1},\cdots, f_{n-1} \rangle$
for $0 <k < n-1,$
\item $d_{0}\gamma =  \langle f_1,\cdots, f_{n-1} \rangle,$
\item $d_{n}\gamma = \langle f_0,f_1,\cdots, f_{n-2} \rangle,$
\item $s_{k}\gamma =  \langle f_0,f_1,\cdots, f_{k-1}, 1_{k}, f_{k}, \cdots,  f_{n-1} \rangle,$
where $1_{k}: A_{k} \longrightarrow A_{k}$ is the identity map on the object $A_{k}$, $f \circ g $ denotes composition of morphisms and $0<k<n-1$.
\item $s_{n}\gamma = \langle f_0,f_1,\cdots, f_{n-1}, 1_{n} \rangle$
\end{enumerate}

\noindent {\bf Theorem.} With the above definitions of face and degeneracy operators, $NerveC$ is a simplicial object.\\

\noindent{\bf Proof.} It suffices to prove the identities
\begin{enumerate}
 \item $d_i d_j = d_{j-1} d_i $ for $i < j.$
 \item $s_i s_j = s_{j+1} s_{i}$ for $i \le j.$
 \item $d_i s_j = s_{j-1} d_i $ for $i < j.$
 \item $d_i s_j = Id$ for $i=j, j+1.$
 \item $d_i s_j = s_j d_{i-1}$ for $i > j+1.$
 \end{enumerate}

\noindent The reader will have no difficulty verifying these identities. We illustrate with one instance of identity $5.$
$$d_{5}s_{3}\gamma = d_{5} \langle f_{0}, f_{1}, f_{2}, 1_{3}, f_{3}, f_{4}, f_{5}, \cdots, f_{n} \rangle$$
$$=  \langle f_{0}, f_{1}, f_{2}, 1_{3}, f_{4} \circ f_{3} , f_{5}, \cdots, f_{n} \rangle,$$
while
$$s_{3}d_{4}\gamma = s_{3} \langle f_{0}, f_{1}, f_{2}, f_{4} \circ f_{3}, f_{5}, \cdots, f_{n} \rangle$$
$$=  \langle f_{0}, f_{1}, f_{2}, 1_{3}, f_{4} \circ f_{3}, f_{5}, \cdots, f_{n} \rangle.$$
Thus we have verified that $d_{5}s_{3} = s_{3}d_{4}.$
We leave the rest of the verifications to the reader. Note that the index $k$ in $1_{k}$ always refers to the domain $A_{k}$ of the corresponding morphism.
This completes the proof. \fbox{}\\

\noindent{\bf Remark.} Note that a sequence of morphisms $ \gamma = \langle f_0,f_1,\cdots, f_{n-1} \rangle,$ as discussed above corresponds to the structure of an $n$-simplex as illustrated in Figure~\ref{morphism}.
In this figure we illustrate how $\gamma = \langle f,g,h \rangle$ corresponds to a tetrahedron, with the edges corresponding to $f,g,h$ and the compositions $hg,gf$ and $(hg)f = h(gf).$ The faces correspond to the faces we have described algebraically above via the face operators $d_{k}.$\\

By taking a functor $F$ from $C$ to a module category $M$ we have the $Nerve(FC)= FNerveC$ and can associate the free abelian simplicial module $ZFNerveC$ to the  simplicial object $FNerveC.$  Such a functor $F$ will be referred to as a {\it presheaf} on the category
$C$ (this is a standard terminology).
We shall now see that the homology or cohomology of this simplicial module
is then a version of the homology or cohomology of the category $C$ with coefficients from the functor $F$ on $C.$ We can also consider directly the homotopy type of $ZFNerveC$,
since it will satisfy the Kan extension condition. We now give the 
details of this construction.\\
  
Let  $F:C \longrightarrow M$ be a functor from $C$ to a module category $M$ such that a terminal object $T$ in $C$ goes to $F(T)$ , a one dimensional module in the category $M.$ Morphisms from $F(T)$ to objects $O$ in $M$ are equivalent to giving an element of the module $O.$ We take $A_{-1} = T$ in the above construction of sequences in the category $C$ (extending the grading one step). Thus $F(NerveC)$ consists in sequences
$ \gamma = \langle g_{-1}; g_0,g_1,\cdots, g_{n-1} \rangle$ of morphisms in  $M$ where $g_k = F(f_k)$ and we can replace $g_{-1}: F(T) \longrightarrow F(A_{-1})$ with an element $\lambda \in F(A_{0}).$
Then we can write 
$$\gamma = \langle \lambda ; g_{0},\cdots, g_{n-1} \rangle,$$
\begin{enumerate}
 \item $d_{0}\gamma =  \langle g_{0}(\lambda) ; g_1,\cdots, g_{n-1} \rangle,$ 
 \item $d_{k}\gamma = \langle \lambda ; g_0,g_1,\cdots g_{k-2}, g_{k} \circ g_{k-1} , \cdots , g_{n-1}  \rangle, $ for $k = 1,\cdots, n-1$
 \item$d_{n}\gamma = \langle \lambda; g_{0},\cdots g_{n-2} \rangle $
 \end{enumerate}
 \noindent and the degeneracies as before.\\
 
\noindent By {\it definition} we identify
$$ \gamma = \langle \lambda ; g_0,g_1,\cdots, g_{n-1} \rangle = \langle \lambda_{0}, \lambda_{1}, \cdots \lambda_{n} \rangle$$ where $$\lambda_{k} = g_{k-1} \circ g_{k-2} \circ \cdots \circ g_{0} (\lambda).$$
Note that $\lambda_{0} = \lambda.$
Thus while a generating simplex for $ZFNerveC$ is a module element followed by  a sequence of composable maps,  we regard two such simplices to be equal if the corresponding sequences of composed module elements are equal, term by term. Note that the $i$-th face operator then becomes just the removal of the $i$-th element in the sequence of $\lambda_{k},$ and degeneracy operators make repetitions just as in an abstract simplex.\\

With this description of the image simplicial object
$F(NerveC),$ we define  $ZFNerveC$ to be the module simplicial object whose $n+1$-th grading is the free module generated by $F(NerveC_{n}).$ Then $ZFNerveC$ is a simplicial module object, satisfying the Kan extension condition and we can associate its homotopy type to the pair $(C, F)$ and we can define the homology and cohomology of the category $C$ with coefficients in $F$ to be the cohomology and homology of
$FNerveC:$  
$$\pi_{*}(C, F) = H_{*}(C, F) = H_{*}(ZFNerveC)= \pi_{*}(ZFNerveC) $$
and
$$H^{*}(C, F) = H^{*}(ZFNerveC).$$ 
We have added homotopy groups to the above formulas since homotopy and homology coincide for 
group complexes. The main point is that there is a single homotopy type behind these groups. Measuring the homology of a category with coefficients in a presheaf $F: C \longrightarrow M$ is preceded by the construction of a simplicial object with a homotopy type that carries this homology. \\

\noindent{\bf The $n$-Simplex Category}. The {\it $n$-Simplex Category $Cat[n]$} is the category generated by a single abstract $n$-simplex $\langle 0,1, \cdots , n \rangle$ and its face maps. The objects of this category
are the given $n$-simplex and all its faces. Each face has face maps to smaller faces on down to the $0$-simplices $\{ \langle 0 \rangle, \langle 1 \rangle , \cdots , \langle n \rangle \}.$  View Figure~\ref{simplexcat} for a depiction of $Cat[2],$ the $2$-Simplex Category. This figure shows how the $2$-Simplex Category can be superimposed on the edges of a geometric $2$-simplex, and how the face morphisms from a given face of the $n$-simplex can be seen as new edges that subdivide that simplex into smaller simplices of the same dimension. Thus in the figure we see that the $2$-dimensional simplex is subdivided into five smaller $2$-simplices and each 
$1$-dimensional face is subdivided into two $1$-dimensional simplices. Just above this diagram in the figure, the reader will see a depiction of the {\it barycentric subdivision} of a geometric $2$-simplex. {\it Each simplex in the barycentric subdivision corresponds to a sequence of morphisms in the simplex category.}  This is a general fact about $Cat[n]$: Its morphism structure encodes the simplicies of the barycentric subdivision of an $n$-simplex. Note that in the Figure~\ref{simplexcat}
that the $2$-simplices $R,S,T,U,V,W$ of the subdivision correspond respectively to the composable morphism sequences 
$$\langle d_{0}, d_{0} \rangle, \langle d_{1}, d_{0} \rangle, \langle d_{1}, d_{1} \rangle, \langle d_{2}, d_{1} \rangle, \langle d_{2}, d_{0} \rangle, \langle d_{0}, d_{1} \rangle .$$
In each case the third side of the corresponding triangle is the composition of the two morphisms in the sequence. The simplicial identities make sure that the triangles are fitted together properly. The upshot of this observation is the following Theorem.\\

\noindent {\bf Subdivision Theorem.} Let $\nabla[n]$ denote the pre-simplicial set generated by an abstract $n$-simplex. ($\nabla[n]$ satisifies all the identities for the face operators $d_{i}$, but has no degeneracies.) 
Let $$F: Cat[n] \longrightarrow M$$ be a functor from the $n$-simplex category to a module category $M.$ Let $Subdiv(\nabla[n])$ denote the pre-simplicial set induced by the barycentric subdivision of
the underlying abstract $n$-simplex $\langle 0,1,2,\cdots , n \rangle$ of $\nabla[n].$ Let $F(Subdiv(\nabla[n]))$ denote the pre-simplicial module object induced by the functor $F.$ Then the simplicial homology of  $F(Subdiv(\nabla[n]))$ is isomorphic to the homology of the category $Cat[n]$ with coefficients in $F.$ That is
$$H_{*}(Cat[n], F) \cong H_{*}(F(Subdiv(\nabla[n])).$$ Furthermore, let $F(\nabla[n])$ denote the pre-simplicial module induced on this pre-simplicial object by the functor $F.$ Then 
$$H_{*}(F(Subdiv(\nabla[n])) \cong H_{*}(F(\nabla[n])).$$ Thus we have an isomorphism of the homology of the subdivision with the homology of the corresponding pre-simplicial module object, and both are isomophic to the homology of the category $Cat[n]$ with coefficients in the functor $F.$\\

\noindent{\bf Proof.} The proof of this theorem is a direct consequence of the correspondence of simplices in the barycentric subdivision with sequences of morphisms in the category $F(Cat[n]).$ Note that we formulated
the definition of $H_{*}(Cat[n], F)$ in terms of an assigment of an element $\lambda$ to each simplex and the boundary operator corresponds exactly to the boundary of the simplex. These points have been illustrated 
in Figure~\ref{morphism} and Figure~\ref{simplexcat}. The second part of the proof is a generalization of the invariance of classical homology under subdivision and has been shown in \cite{Jozef1,Jozef2,Wang}.
This completes the proof.//\\

\noindent {\bf Remark.} In stating this Theorem we have been careful to emphasize the pre-simplicial objects to which it refers. However there is a full simplicial object and associated simplicial module that is  relevant to this
result. In this section we have described the simplicial object $NerveC$ where $C$ is a small category. Here we have $NerveCat[n]$ and the simplicial module object $F(NerveCat[n])$ generated by sequences
$\langle \lambda; g_0,g_1,\cdots,g_n \rangle$ where $\lambda$ is an element of the domain of $g_0$ and the $g_k$ are module morphisms. We have, by definition, that $$H_{*}(Cat[n], F) \cong H_{*}(F(NerveCat[n])).$$
Thus we can identify the homotopy type of $ZF(NerveCat[n])$ in back of this Theorem.\\

\noindent {\bf Remark.} The Subdivision Theorem generalizes to {\it simplicial categories} by which we mean categories that have the structure of a simplicial complex in the same sense that a simplex category has the structure of a simplex.
We will study this generalization in a sequel to the present paper. If we were given an extension of the functor $F$ to the simplicial object $\Delta[n]$ we could say more.
We could then work with a simplex category $F(Cat[n])$ over a module category and an associated simplicial module object $F(\Delta[n]).$ We have that the simplicial object associated with the nerve of the category (represented by sequences $\langle \lambda ; g_1 , g_2, \cdots g_n \rangle$ as above) corresponds to the barycentric subdivision of  the simplicial module object $F(\Delta[n])$ and is identical with $F(Subdiv(\Delta[n])).$ In fact,  
$ZF(\Delta[n])$ and $ZF(Subdiv(\Delta[n]))$ have the same homotopy type, from which one can also conclude the above Theorem. In a sequel to this paper, we will reformulate the Theorem in this way.  In the next section, we shall apply the ideas of this  construction to Khovanov homology by showing how it is isomorphic with a categorical homology for a category that is specifically associated with a knot or link diagram.\\

\section{Khovanov Homology and the Cube Category}
Khovanov homology for a knot or link diagram $K$ is constructed from a cube category constructed from the Kauffman bracket states \cite{K,K1,BN1,BN2} of $K.$
The objects in this cube category are the states of the bracket as shown in Figure~\ref{Figure 2}. To each state circle is associated a module $V$ and to each state is associated a module $V(S)$ that is the 
tensor product with one $V$ for each state loop. This assigment gives a functor from an $n$-cube category to a category of modules. We shall futher describe Khovanov homology by showing how that cube category embeds in a simplicial category. Then it can be explained how the methods described here can be applied.\\

 Examine Figure~\ref{Figure 2} and Figure~\ref{Figure 3}. In Figure~\ref{Figure 2} we show all the 
standard bracket states for the trefoil knot with arrows between them whenever the state at the output of the arrow is obtained from the state at the input of the arrow by a single smoothing of a site of type $A$ to a site of type $B$. In Figure~\ref{Figure 3} we illustrate the  cube  category (the states are arranged in the form of a cube)
by replacing the states in Figure~\ref{Figure 2} by symbols $\langle XYZ \rangle$ where each literal is either $A$ or $B.$
A typical generating morphism in the $3$-cube category is $$\langle ABA\rangle \longrightarrow \langle BBA\rangle .$$
\bigbreak

Let $\mathcal{S}(K)$ denote a category 
associated with the states of the bracket for a diagram $K$ whose objects are the states, with sites
labeled $A$ and $B$ as in Figure~\ref{Figure 2}. A (generating) morphism in this category is an arrow from a state with a given
number of $A$'s to a state with {\it one less} $A.$ We also have identity maps from each object to itself, but these are not indicated in the figures.
Note that the objects in this category are {\it states}; they are collections of circles with site labels $A$ and $B.$
\bigbreak

 \noindent Let ${\mathcal D}^{n} = \{A,B \}^{n}$ be the $n$-cube category whose objects are the $n$-sequences from
the set $\{A, B \}$ and whose non-identity morphisms are arrows from sequences with greater numbers of $A$'s
to sequences with fewer numbers of $A$'s. Thus ${\mathcal D}^{n}$ is equivalent to the poset category
of subsets of $\{ 1,2,\cdots n \}.$ We make a functor $\mathcal{R}: {\mathcal D}^{n} \longrightarrow \mathcal{S}(K)$
for a diagram $K$ with $n$ crossings as follows.
We map objects in the cube category  to bracket states
by choosing to label the crossings of the diagram $K$ from the set  $\{ 1,2,\cdots n \},$ and letting
this  functor take abstract $A$'s and $B$'s in the cube category to smoothings at those crossings of type $A$ or type $B.$ Thus each object in the cube category is associated with a unique state of $K$
when $K$ has $n$ crossings. By the same token, we define a functor 
$\mathcal{S}:\mathcal{S}(K) \longrightarrow {\mathcal D}^{n} $ by associating an object (seqeuence of $A$'s and $B$'s to each state
and morphisms between sequences corresponding to the state smoothings. These are inverse functors.
\bigbreak

In Figure~\ref{cubeplex}, Figure~\ref{cubecat} and Figure~\ref{simplexcat} we show how one can embed the cube category in a simplex category. The combinatorics of this embedding is as follows. Let 
$$\langle *AAA \cdots A \rangle = \langle 012\cdots n\rangle$$ stand for an abstract
$n$-simplex, where there are $n$ copies of $A.$ Define face maps as follows (here indicated for $n = 3$).
$$d_{0}\langle *AAA \rangle = \langle \hat{0}123 \rangle = \langle \sharp AAA \rangle,$$
$$d_{1}\langle *AAA \rangle = \langle 0\hat{1}23 \rangle = \langle * BAA \rangle,$$
$$d_{2}\langle *AAA \rangle = \langle 01\hat{2}3 \rangle = \langle * ABA \rangle,$$
$$d_{3}\langle *AAA \rangle = \langle 012\hat{3} \rangle = \langle * AAB \rangle.$$
Each $B$ is regarded as an eliminated sites in the abstract simplex.  For example,
$$d_{3}\langle *BBA \rangle = \langle *BBB \rangle.$$ As Figure~\ref{cubeplex} and Figure~\ref{cubecat} show, this gives a natural embedding of the cube category into a simplex category.\\

Let ${\mathcal M }$ be a module category and let $V$ be a chosen module. We shall associate a {\it state module} $V(S)$ to each state $S$ of the knot $K$ by taking the tensor product of one copy of $V$ for each loop in $S.$ If $V$ is equipped with maps $m: V \otimes V \longrightarrow V$ and $\Delta: V \longrightarrow V \otimes V,$ then we can use them to define a map of modules $V(S) \longrightarrow V(S')$ whenever there is a morphism $f:S \longrightarrow S' $ in the category $\mathcal{S}(K).$ To see this, note that a generating morphism consists in the re-smoothing of a state as illustrated in Figure~\ref{Figure 2}, and that such a re-smoothing acts either on two circles to form one circle, or on one circle to form two circles. If we apply the maps $m$ and $\Delta$ correspondingly, we obtain a well-defined map on the corresponding state modules. The result then follows by composition of morphisms. We want the resulting morphisms between state modules to be independent of the factorization of a morphism into generating morphisms. This leads to algebraic conditions on the maps
$m$ and $\Delta.$ They must define a {\it Frobenius algebra} structure on the module $V.$ See \cite{K,K1} for the interesting details of this construction. A relevant example of such a Frobenius algebra is 
$V = k[x]/(x^{2})$ where $k$ is a ground ring with unit and $m(x \otimes 1) = m(1 \otimes x) = x, m(x \otimes x ) = 0, \Delta(x) = x \otimes x, \Delta(1) = 1 \otimes x + x \otimes 1.$ Given a Frobenius algebra $V$ we have 
a well-defined functor $\hat{V}: \mathcal{S}(K) \longrightarrow \mathcal{M},$ and with that we have the functor $$F: {\mathcal D}^{n} = \{A,B \}^{n} \longrightarrow {\mathcal M }$$ obtained by composing
$\mathcal{R}$ with $\hat{V}.$\\

 Khovanov homology is usually defined directly in terms of the cubical description, but we will give the definition here in terms of the simplicial embedding.\\

For the functor $F$ we first construct a semisimplicial object ${\mathcal C(F)}$ over ${\mathcal M }$,
where a semisimplicial object is a simplicial object without degeneracies. For $k \ge 0$
we set $${\mathcal C(F)}_{k} = \oplus_{v  \in {\mathcal D}^{n}_{k}} {\mathcal F}(v)$$
where ${\mathcal D}^{n}_{k}$
denotes those objects in the cube category with $k$ $A$'s. Note that we are indexing dually
to the upper indexing in the Khovanov homology sections of this paper where we counted the number of
$B$'s in the states.
\bigbreak

We introduce face operators (partial boundaries in our previous terminology)
$$d_{i} : {\mathcal C(F)}_{k} \longrightarrow {\mathcal C(F)}_{k-1}$$ for $0 \le i \le k$ with $k \ge 1$ as
follows: $d_{i}$ is trivial for $i = 0$ and otherwise $d_{i}$ acts on ${\mathcal F}(v)$ by the map
${\mathcal F}(v) \longrightarrow {\mathcal F}(v')$ where $v'$ is the object resulting from replacing
the $i$-th $A$ by $B.$ The operators $d_{i}$ satisfy the usual face relations of simplicial theory:
$$d_{i} d_{j} = d_{j-1} d_{i}$$ for $i <  j.$ Note that this algebraic description is exactly parallel to our description of the embedding of the cube category in a simplex category.
\bigbreak

We now expand ${\mathcal C}(F)$ to a simplicial object ${\mathcal S}(F)$ over ${\mathcal M }$ by applying freely
degeneracies to the ${\mathcal F}(v)$'s. Thus
$${\mathcal S}(F)_{m} = \oplus_{v \in {\mathcal D}^{n}_{k}, k + t = m} \,\, s_{i_{1}} \cdots  s_{i_{t}} {\mathcal F}(v)$$
where $ m  >  i_{1}  >  \cdots  > i_{t} \ge 0 $ and these degeneracy operators are applied freely modulo the usual (axiomatic) relations among themselves and with the face operators. Then ${\mathcal S}(F)$
has degeneracies via formal application of degeneracy operators to these forms, and has face operators
extending those of ${\mathcal C}({\mathcal F}).$ \\

It is at this point we should remark that in our knot theoretic construction
there is only at this point an opportunity for formal extension of degeneracy operators above the number of crossings in the given knot or link diagram since to make specific degeneracies would involve the
creation of new diagrammatic sites. Nevertheless, we have created a simplicial abelian group that sits in back of Khovanov homology. The simplcial object ${\mathcal S}(F)$ is a complex with a homotopy type satisfying the hypotheses of the Dold-Kan correspondence. Thus we know that the homotopy type of ${\mathcal S}(F)$ corresponds to the chain-homotopy type of $N[{\mathcal S}(F)]$ where
$N$ denotes the associated normalized chain complex. We also know that the homology of $N[{\mathcal S}(F)]$ is isomorphic with the Khovanov homology of our original knot. The next paragraphs detail this correspondence.\\

When the functor $F: {\mathcal D}^{n} \longrightarrow {\mathcal M}$ goes to an abelian category ${\mathcal M},$
as in our knot theoretic case, we can recover the homology groups via
$$H_{\star} N {\mathcal S}(F) \cong H_{\star} {\mathcal C}(F)$$
where $N {\mathcal S}({\mathcal F})$ is the normalized chain complex of ${\mathcal S}({\mathcal F}).$ This completes
the abstract simplicial description of this homology.
\bigbreak

\noindent {\bf Theorem.} Khovanov homology of the knot or link $K$ is given by $H_{\star} {\mathcal C}(F)$ where $F$ is the Frobenius functor on the embedding of the Khovanov cube category in a simplex category.\\

\noindent {\bf Proof.} The Theorem is a tautology based on the original definition of Khovanov homology using the cube category. The face maps in the simplicial description correspond exactly to the partial boundaries obtained by re-smoothing crossings in cubical depiction. //\\

With this description in hand, we can apply the Dold-Kan Theorem. We have that $KN {\mathcal S}(F)$ is a simplicial object such that chain homotopies of $N {\mathcal S}({\mathcal F})$
correspond to  homotopies of $KN {\mathcal S}(F).$ We also know, by the Dold-Kan correspondence that 
 $$KN {\mathcal S}(F) =  {\mathcal S}(F).$$ Thus we know that the simplicial object ${\mathcal S}(F)$ carries the homotopy properties relevant to Khovanov homology.
 {\it The chain homotopy type of $N{\mathcal S}(F)$  is equivalent to the homotopy type of ${\mathcal S}(F),$ and this chain homotopy type contains all the information about Khovanov homology. So we have found a homotopy interpretation for Khovanov homology in terms of its own cube category with respect to a Frobenius functor $ F.$}
 
 \subsection{Khovanov Homology,  Categorical Homology and Homotopy}
 We have the defining functor $F$, for Khovanov homology, from the cube category to a Frobenius module category. This extends to a functor $F$ defined on a simplex category $Cat[n].$ We have defined Khovanov 
 homology via $KhoH_{*}(K)= H_{*}(F(\nabla[n])$ in the terminology of the Subdivision Theorem of Section 6. By this Theorem we have that  $H_{*}(F(\nabla[n])) \cong H_{*}(F(NerveCat[n]))$ where $F(NerveCat[n])$ is the 
 simplicial module generated by sequences $\langle \lambda; g_0,g_1,\cdots,g_n \rangle$ where $\lambda$ is in the domain of $g_0$ and the $g_{k}$ are sequences of composable morphisms in the module category, induced by the corresponding maps in the Khovanov complex for the link $K.$
We can take the homotopy type of the simplicial module object $F(NerveCat[n])$ and regard it as the appropriate homotopy structure behind Khovanov homology. \\

We end here at the beginning of an investigation into the homotopy structure of Khovanov homology. The simplical module object $F(NerveCat[n])$ is a good place to begin the homotopy theory, and that will be the subject for the next paper.\\

\noindent {\bf Questions and Directions.} Here are a number of questions and remarks that arise directly from the present paper.
\begin{enumerate}
\item According to \cite{May} (Theorem 24.5, page 106): If G is a connected Abelian group complex, and $\pi_{n} = \pi_{n}(G)$, then G has the homotopy type of the (infinite) Cartesian product of Eilenberg-MacLane spaces
$\Pi_{n=1}^{\infty} K(\pi_{n}, n).$ This means that all homotopy types discussed in this paper are products of Eilenberg-MacLane spaces. The constructions of Lifshitz and Sarkar \cite{LS1,LS2,LS3} are not simply products, and this makes it even more interesting to see the relationship between those constructions and the simplicial constructions of the present paper.
\item Are there smaller and more accessible spatial models than the simplicial models indicated here? In particular, let $K$ be a link diagram. Let $\Lambda(K)$ denote the nerve of the image of the Khovanov category for
$K$ under the Frobenius algebra functor that associates a module to each object of the Khovanov category for $K.$ The homology of $\Lambda(K)$ can be taken as the Khovanov homology of $K.$ The simplicial object
$\Lambda(K)$ is simpler than the free abelian group object $Z\Lambda(K)$ whose homotopy type is a product of Eilenberg-Maclane spaces for the homology of $\Lambda(K).$ In this paper, we have promoted 
$Z\Lambda(K)$ as a possibly useful homotopy type for Khovanov homology of $K.$ But $\Lambda(K)$ already carries this homology. It is a very interesting question to ask about the homotopy type of the realization
$|\Lambda(K)|.$ This is a space that carries the Khovanov homology. What are its homotopy groups? How do they behave under isotopies of the link $K?$ We are indebted to Stephan Klaus for this line of questioning.

\item Is there a more direct construction of a homotopy type for Khovanov homology that depends directly on the embedding of the knot or link in three dimensional space, circumventing the dependence on the structure of 
the link diagram?
\end{enumerate}
 
\section{ Appendix -- Bracket Polynomial, Jones Polynomial and Grading in Khovanov Homology}
The bracket polynomial \cite{KaB} model for the Jones polynomial \cite{JO,JO1,JO2,Witten} is usually described by the expansion
$$\langle \Across \rangle=A \langle \Asmooth \rangle + A^{-1}\langle \Bsmooth \rangle$$
Here the small diagrams indicate parts of otherwise identical larger knot or link diagrams. The two types of smoothing (local diagram with no crossing) in
this formula are said to be of type $A$ ($A$ above) and type $B$ ($A^{-1}$ above).

$$\langle \bigcirc \rangle = -A^{2} -A^{-2}$$
$$\langle K \, \bigcirc \rangle=(-A^{2} -A^{-2}) \langle K \rangle $$
$$\langle \Rcurl \rangle=(-A^{3}) \langle \Arc \rangle $$
$$\langle \Lcurl \rangle=(-A^{-3}) \langle \Arc \rangle $$
One uses these equations to normalize the invariant and make a model of the Jones polynomial.
In the normalized version we define $$f_{K}(A) = (-A^{3})^{-wr(K)} \langle K \rangle / \langle \bigcirc \rangle $$
where the writhe $wr(K)$ is the sum of the oriented crossing signs for a choice of orientation of the link $K.$ Since we shall not use oriented links
in this paper, we refer the reader to \cite{KaB} for the details about the writhe. One then has that $f_{K}(A)$ is invariant under the Reidemeister moves
(again see \cite{KaB}) and the original Jones polynomial $V_{K}(t)$ is given by the formula $$V_{K}(t) = f_{K}(t^{-1/4}).$$ The Jones polynomial
has been of great interest since its discovery in 1983 due to its relationships with statistical mechanics, due to its ability to often detect the
difference between a knot and its mirror image and due to the many open problems and relationships of this invariant with other aspects of low
dimensional topology.
\bigbreak

\noindent {\bf The State Summation.} In order to obtain a closed formula for the bracket, we now describe it as a state summation.
Let $K$ be any unoriented link diagram. Define a {\em state}, $S$, of $K$  to be the collection of planar loops resulting from  a choice of
smoothing for each  crossing of $K.$ There are two choices ($A$ and $B$) for smoothing a given  crossing, and
thus there are $2^{c(K)}$ states of a diagram with $c(K)$ crossings.
In a state we label each smoothing with $A$ or $A^{-1}$ according to the convention
indicated by the expansion formula for the bracket. These labels are the  {\em vertex weights} of the state.
There are two evaluations related to a state. The first is the product of the vertex weights,
denoted $\langle K|S \rangle .$
The second is the number of loops in the state $S$, denoted  $||S||.$

\noindent Define the {\em state summation}, $\langle K \rangle $, by the formula

$$\langle K \rangle  \, = \sum_{S} <K|S> \delta^{||S||}$$
where $\delta = -A^{2} - A^{-2}.$
This is the state expansion of the bracket. It is possible to rewrite this expansion in other ways. For our purposes in
this paper it is more convenient to think of the loop evaluation as a sum of {\it two} loop evaluations, one giving $-A^{2}$ and one giving
$-A^{-2}.$ This can be accomplished by letting each state curve carry an extra label of $+1$ or $-1.$ We describe these {\it enhanced states} below. But before we do this, it will be useful for the reader to examine Figure~\ref{Figure 2}. In Figure~\ref{Figure 2} we show all the states for the right-handed trefoil knot, labelling the sites
with $A$ or $B$ where $B$ denotes a smoothing that would receive $A^{-1}$ in the state expansion.
\bigbreak

Note that in the state enumeration in Figure~\ref{Figure 2} we have organized the states in tiers so that the state
that has only $A$-smoothings is at the top and the state that has only $B$-smoothings is at the bottom.
\bigbreak

\begin{figure}
     \begin{center}
     \begin{tabular}{c}
     \includegraphics[width=6cm]{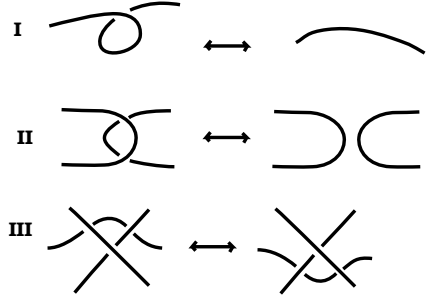}
     \end{tabular}
     \caption{\bf Reidemeister Moves}
     \label{Figure 1}
\end{center}
\end{figure}

\noindent {\bf Changing Variables.} Letting $c(K)$ denote the number of crossings in the diagram $K,$ if we replace $\langle K
\rangle$ by
$A^{-c(K)} \langle K \rangle,$ and then replace $A^2$ by $-q^{-1},$ the bracket is then rewritten in the
following form:
$$\langle \Across \rangle=\langle \Asmooth \rangle-q\langle \Bsmooth \rangle $$
with $\langle \bigcirc\rangle=(q+q^{-1})$.
It is useful to use this form of the bracket state sum
for the sake of the grading in the Khovanov homology (to be described below). We shall
continue to refer to the smoothings labeled $q$ (or $A^{-1}$ in the
original bracket formulation) as {\it $B$-smoothings}.
\bigbreak

We catalog here the resulting behaviour of this modified bracket under the Reidemeister moves.
$$\langle \bigcirc \rangle = q + q^{-1}$$
$$\langle K \, \bigcirc \rangle=(q + q^{-1}) \langle K \rangle $$
$$\langle \Rcurl \rangle=q^{-1} \langle \Arc \rangle $$
$$\langle \Lcurl \rangle= -q^{2} \langle \Arc \rangle $$
$$\langle \raisebox{-0.25\height}{\includegraphics[width=0.5cm]{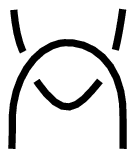}}\rangle =
-q \langle  \raisebox{-0.50\height}{\includegraphics[width=0.5cm]{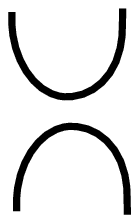}}\rangle$$
$$\langle \raisebox{-0.25\height}{\includegraphics[width=0.8cm]{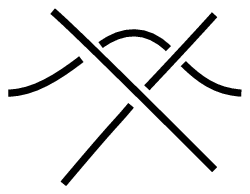}}\rangle =
\langle  \raisebox{-0.25\height}{\includegraphics[width=0.8cm]{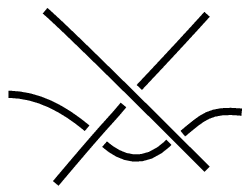}}\rangle$$

It follows that if we define $$J_{K} = (-1)^{n_{-}} q^{n_{+} - 2n_{-}} \langle K \rangle,$$
where $n_{-}$ denotes the number of negative crossings in $K$ and $n_{+}$ denotes the number
of positive crossings in $K$, then $J_{K}$ is invariant under all three Reidemeister moves.
Thus $J_{K}$ is a version of the Jones polynomial taking the value $q + q^{-1}$ on an unknotted circle.
\bigbreak

\noindent {\bf Using Enhanced States.}
We now use the convention of {\it enhanced
states} where an enhanced state has a label of $1$ or $-1$ on each of
its component loops. We then regard the value of the loop $q + q^{-1}$ as
the sum of the value of a circle labeled with a $1$ (the value is
$q$) added to the value of a circle labeled with an $-1$ (the value
is $q^{-1}).$ We could have chosen the less neutral labels of $+1$ and $X$ so that
$$q^{+1} \Longleftrightarrow +1 \Longleftrightarrow 1$$
and
$$q^{-1} \Longleftrightarrow -1 \Longleftrightarrow x,$$
since an algebra involving $1$ and $x$ naturally appears later in relation to Khovanov homology. It does no harm to take this form of labeling from the
beginning. The use of enhanced states for formulating Khovanov homology was pointed out by Oleg Viro in
\cite{Viro}.
\bigbreak

\begin{figure}
     \begin{center}
     \begin{tabular}{c}
     \includegraphics[width=7cm]{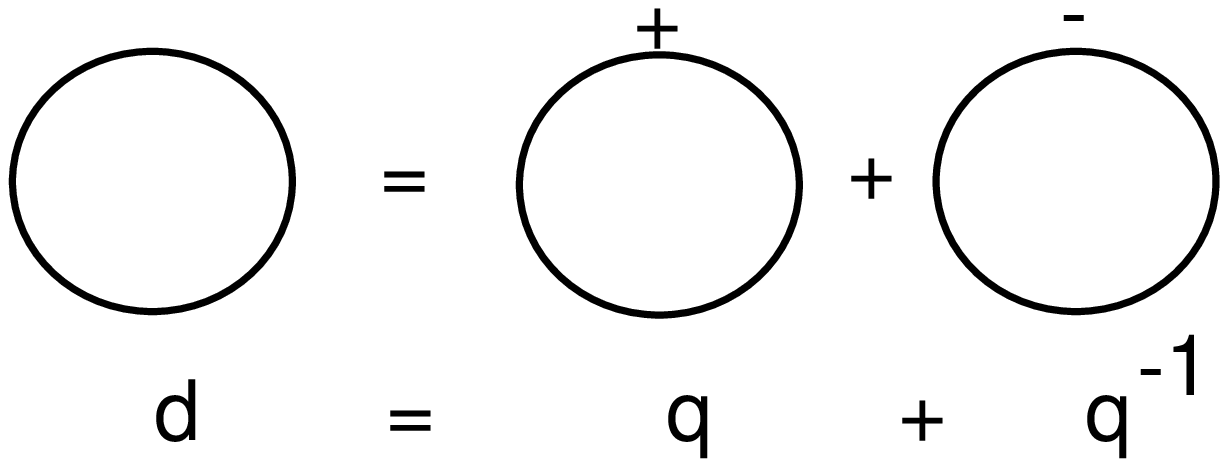}
     \end{tabular}
     \caption{\bf Enhanced States for One Loop}
    \label{enhancedloop}
\end{center}
\end{figure}

Consider the form of the expansion of this version of the
bracket polynonmial in enhanced states. We have the formula as a sum over enhanced states $s:$
$$\langle K \rangle = \sum_{s} (-1)^{i(s)} q^{j(s)} $$
where $i(s)$ is the number of $B$-type smoothings in $s$ and $j(s) = i(s) + \lambda(s)$, with $\lambda(s)$ the number of loops  labeled $1$ minus
the number of loops labeled $-1$ in the enhanced state $s.$
\bigbreak

Two key motivating ideas are involved in finding the Khovanov invariant. First of all, one would like to {\it categorify} a link polynomial such as
$\langle K \rangle.$ There are many meanings to the term categorify, but here the quest is to find a way to express the link polynomial
as a {\it graded Euler characteristic} $\langle K \rangle = \chi_{q} \langle {\mathcal H}(K) \rangle$ for some homology theory associated with $\langle K \rangle.$
\bigbreak

To see how the Khovanov grading arises, consider the form of the expansion of this version of the
bracket polynomial in enhanced states. We have the formula as a sum over enhanced states $s:$
$$\langle K \rangle = \sum_{s} (-1)^{i(s)} q^{j(s)} $$
where $i(s)$ is the number of $B$-type smoothings in $s$, $\lambda(s)$ is the number of loops in $s$ labeled $1$ minus the number of loops
labeled $X,$ and $j(s) = i(s) + \lambda(s)$.
This can be rewritten in the following form:
$$\langle K \rangle  =  \sum_{i \,,j} (-1)^{i} q^{j} dim({\mathcal C}^{ij}) $$
where we define ${\mathcal C}^{ij}$ to be the linear span (over the complex numbers or over the integers or the integers modulo
two for other contexts) of the set of enhanced states with
$i(s) = i$ and $j(s) = j.$ Then the number of such states is the dimension $dim({\mathcal C}^{ij}).$
\bigbreak

\noindent We would like to have a  bigraded complex composed of the ${\mathcal C}^{ij}$ with a
differential
$$\partial:{\mathcal C}^{ij} \longrightarrow {\mathcal C}^{i+1 \, j}.$$
The differential should increase the {\it homological grading} $i$ by $1$ and preserve the
{\it quantum grading} $j.$
Then we could write
$$\langle K \rangle = \sum_{j} q^{j} \sum_{i} (-1)^{i} dim({\mathcal C}^{ij}) = \sum_{j} q^{j} \chi({\mathcal C}^{\bullet \, j}),$$
where $\chi({\mathcal C}^{\bullet \, j})$ is the Euler characteristic of the subcomplex ${\mathcal C}^{\bullet \, j}$ for a fixed value of $j.$
\bigbreak

\noindent This formula would constitute a categorification of the bracket polynomial. In fact,  {\it the original Khovanov differential $\partial$ is uniquely determined by the restriction that $j(\partial s) = j(s)$ for each enhanced state $s$.}  In particular one takes the multiplication induced by the algebra ${\cal A} = k[x]/(x^2)$ so that $m(1\otimes x) = m(x \otimes 1) = x, m(1 \otimes 1) = 1, m(x \otimes x) = 0,$ and the comultiplication
$\Delta: {\cal A} \longrightarrow {\cal A} \otimes {\cal A},$ with $\Delta(1) = 1 \otimes x + x \otimes 1, \Delta(x) = x \otimes x.$ These operations of multiplcation and comutiplication act on single or double loop enhanced states
labeled according to the algebra. The enhanced states generate the chain complex and the local differentials are the possible multiplications or comultiplications that change one $A$-type smoothing to a $B$-type smoothing.
The reader can start here and translate these differerentials into the simplicial structures we have discussed earlier in the paper.\\

\noindent Since $j$ is
preserved by the differential, these subcomplexes ${\mathcal C}^{\bullet \, j}$ have their own Euler characteristics and homology. We have
$$\chi(H({\mathcal C}^{\bullet \, j})) = \chi({\mathcal C}^{\bullet \, j}) $$ where $H({\mathcal C}^{\bullet \, j})$ denotes the homology of the complex
${\mathcal C}^{\bullet \, j}$. We can write
$$\langle K \rangle = \sum_{j} q^{j} \chi(H({\mathcal C}^{\bullet \, j})).$$
The last formula expresses the bracket polynomial as a {\it graded Euler characteristic} of a homology theory associated with the enhanced states
of the bracket state summation. This is the categorification of the bracket polynomial. Khovanov proves that this homology theory is an invariant
of knots and links (via the Reidemeister moves of Figure~\ref{Figure 1}), creating a new and stronger invariant than the original Jones polynomial.
\bigbreak

\noindent{\bf Remark on Grading and Invariance.} We showed how the 
bracket, using the variable $q$, behaves under Reidemeister moves. These formulas correspond to how 
the invariance of the homology works in relation to the moves. We have that 
$$J_{K} = (-1)^{n_{-}} q^{n_{+} - 2n_{-}} \langle K \rangle,$$
where $n_{-}$ denotes the number of negative crossings in $K$ and $n_{+}$ denotes the number
of positive crossings in $K.$ $J(K)$ is invariant under all three Reidemeister moves. The corresponding formulas for Khonavov homology are as follows
$$J_{K} = (-1)^{n_{-}} q^{n_{+} - 2n_{-}} \langle K \rangle = $$
$$(-1)^{n_{-}} q^{n_{+} - 2n_{-}} \Sigma_{i,j} (-1)^{i}a^{j} dim(H^{i,j}(K) = $$
$$ \Sigma_{i,j} (-1) ^{i + n_{+} }q^{j + n_{+} - 2n_{-1} }dim(H^{i,j}(K)) = $$
$$\Sigma_{i,j} (-1)^{i} q^{j} dim(H^{i - n_{-}, j - n_{+} + 2n_{-}}(K)).$$
It is often more convenient to define the {\em Poincar\'e polynomial} for Khovanov homology via
$$P_{K}(t, q)  
= \Sigma_{i,j} t^{i} q^{j} dim(H^{i - n_{-}, j - n_{+} + 2n_{-}}(K)).$$ 
The Poincar\'e polynomial is a two-variable polynomial invariant of knots and links, generalizing the Jones polynomial. Each coefficient $$ dim(H^{i - n_{-}, j - n_{+} + 2n_{-}}(K))$$ is an 
invariant of the knot, invariant under all three Reidemeister moves. In fact, the homology groups 
$$H^{i - n_{-}, j - n_{+} + 2n_{-}}(K)$$ are knot invariants. The grading compensations show how the grading of the homology can change from diagram to diagram for diagrams that represent the same knot.
\bigbreak

\end{document}